\documentclass[11pt]{article}
\usepackage{amsmath, amssymb, epsfig, eepic}

\numberwithin{equation}{section}

\newcommand{\eps}{\varepsilon}

\newcommand{\qed}{\hfill $\Box$}

\newcommand{\prob}{\mathbb P}
\newcommand{\expec}{\mathbb E}
\newcommand{\smallsup}[1] {{\scriptscriptstyle{({#1}})}}
\newtheorem{theorem}{Theorem}[section]
\newtheorem{lemma}[theorem]{Lemma}

\newcommand{\eq}{\begin{equation}}
\newcommand{\en}{\end{equation}}
\def\eqalign#1\enalign{
    \begin{align}#1\end{align}
    }
\newcommand{\sss}   { \scriptscriptstyle }
\newcommand{\nn}   { \nonumber}

\newcommand{\sD}{{\scriptscriptstyle D}}

\def\NN{{\mathbb N}}

\def\indic{{\bf 1}}

\newcommand*{\sumu}{\displaystyle\sum}
\newcommand*{\produ}{\displaystyle\prod}
\def\AeN{A_{\varepsilon\kern-.05em,\kern .1em \sN} }
\def\BeN{ B_{\varepsilon\kern-.05em,\kern .1em \sN}}
\def\CeN{C_{\varepsilon\kern-.05em,\kern .1em \sN} }
\def\DeN{ D_{\varepsilon\kern-.05em,\kern .1em \sN} }
\def\EeN{ E_{\varepsilon\kern-.05em,\kern .1em \sN} }
\def\GeN{G_{\varepsilon\kern-.05em,\kern .1em \sN}}

\def\prob{{\mathbb P}}
\def\expec{{\mathbb E}}
\newcommand{\sN}{\sss N}

\title{Distances in random graphs with infinite mean degrees}
\author{Remco van der Hofstad\footnote{Department of Mathematics and
Computer Science, Eindhoven University of Technology, P.O.\ Box
513, 5600 MB Eindhoven, The Netherlands. E-mail: {\tt
rhofstad@win.tue.nl}}\\
Gerard Hooghiemstra \footnote{Delft University of Technology,
Electrical Engineering, Mathematics and Computer Science, P.O. Box 5031, 2600 GA Delft,
The Netherlands. E-mail: {\tt G.Hooghiemstra@ewi.tudelft.nl}} and
Dmitri Znamenski\footnote{EURANDOM, P.O.\ Box 513, 5600 MB
Eindhoven, The Netherlands. E-mail: {\tt znamenski@eurandom.nl}} }

\begin{document}
\date{June 2, 2004}
\maketitle

\begin{abstract}
We study random graphs
with an i.i.d.\  degree sequence of which the tail of the
distribution function $F$ is regularly varying with exponent
$\tau\in (1,2)$. Thus, the degrees have infinite mean.
Such random graphs can serve as models for complex
networks where degree power laws are observed.

The minimal number of edges between two arbitrary nodes, also called the
graph distance or the hopcount, in a graph with $N$ nodes is
investigated when $N\rightarrow \infty$. The paper is part of a sequel
of three papers. The other two papers study the case where
$\tau \in (2,3)$, and $\tau \in (3,\infty),$ respectively.

The main result of this paper is that the graph distance
converges for $\tau\in (1,2)$ to a limit random variable
with probability mass exclusively on the points $2$ and $3$.
We also consider the case where we condition the degrees to be at most
$N^{\alpha}$ for some $\alpha>0.$ For $\tau^{-1}<\alpha<(\tau-1)^{-1}$,
the hopcount converges to 3 in probability, while for
$\alpha>(\tau-1)^{-1}$, the hopcount converges to the same limit as for
the unconditioned degrees. Our results give convincing asymptotics
for the hopcount when the mean degree is infinite, using extreme value theory.
\end{abstract}

\smallskip


\smallskip

\noindent\textit{AMS} 2000 \textit{Subject classifications: Primary 60G70,
Secondary 05C80.} \smallskip


\section{Introduction}
The study of complex networks has attracted considerable attention
in the past decade. There are numerous examples of complex
networks, such as social relations, the World-Wide Web and Internet,
co-authorship and citation networks of scientists, etc.
The topological structure of networks affects the performance
in those networks. For instance, the topology of social networks
affects the spread of information and disease (see e.g., \cite{Strogatz_Nature01}),
while the performance of traffic in Internet depends heavily on
the topology of the Internet.

Measurements on complex networks have shown that many
real networks have similar properties. A first example of such
a fundamental network property is the fact
that typical distances between nodes are small. This is called the
`small world' phenomenon, see the pioneering work of Watts
\cite{Watts} and the references therein.
In Internet, for example, e-mail messages cannot use more than a
threshold of physical links, and if the distances in Internet
would be large, e-mail service would simply break down. Thus, the
graph of the Internet has evolved in such a way that typical
distances are relatively small, even though the Internet is rather
large.

A second, maybe more surprising, property of many networks
is that the number of nodes with degree $k$ falls of as an inverse
power of $k$. This is called a `power law degree sequence'.
In Internet, the power law degree sequence was first observed in
\cite{FFF99}. The observation that many real networks have the above
properties have incited a burst of activity in network modelling.
Most of the models use random graphs as a way to model the
uncertainty and the lack of regularity in real networks.
See \cite{AB02, Newm03} and the references therein for an introduction to
complex networks and many examples where the above two properties
hold.

The current paper presents a rigorous derivation for the random
fluctuations of the  distance between two arbitrary nodes
(also called the geodesic) in a graph with i.i.d.\ degrees with
{\it infinite} mean.  The model with i.i.d.\ degrees is a variation of
the {\it configuration model}, which was originally proposed by
Newman, Strogatz and Watts \cite{NSW00}, where the degrees originate from a given
deterministic sequence. The observed power exponents
are in the range from $\tau=1.5$ to $\tau=3.2$
(see \cite[Table II]{AB02} or \cite[Table II]{Newm03}).
In a previous paper with Van Mieghem \cite{HHV03}, we have
investigated the case $\tau>3$. The case $\tau\in(2,3)$ was
studied in \cite{HHZ04b}. Here we focus on the case
$\tau \in (1,2)$, and study the typical distances between arbitrary
connected nodes. In a forthcoming paper \cite{HHZ04c}, we will survey the
results from the different cases for $\tau$, and investigate the connected
components of the random graphs.

This section is organized as follows. In Section \ref{sec-mod} we
start by introducing the model, and in Section \ref{sec-res} we
state our main results. In Section \ref{sec-heuristics} we explain heuristically the
results obtained. Finally we describe related work in Section \ref{sec-related}.


\subsection{The model}
\label{sec-mod}
Consider an i.i.d.\ sequence $D_1,D_2,\ldots,D_{\sN}$. Assume that
$L_{\sN}=\sum_{j=1}^N D_j$ is even. When $L_{\sN}$ is odd, then we
increase the number of stubs of $D_N$ by $1$, i.e., we replace $D_N$ by
$D_{N}+1$. This change will make hardly any difference in what
follows, and we will ignore it in the sequel.

We will construct a graph in which node $j$ has degree $D_j$ for
all $1\leq j \leq N$. We will later specify the distribution of
$D_j$. We start with $N$ separate nodes and incident to node $j$, we have
$D_j$ stubs which still need to be connected to build the graph.

The stubs are numbered in an arbitrary order from $1$ to $L_{\sN}$. We
now continue by matching at random the first stub with one of the
$L_{\sN}-1$ remaining stubs. Once paired, two stubs form an edge
of the graph. Hence, a stub can be seen as the left or the right
half of an edge. We continue the procedure of randomly choosing
and pairing the next stub and so on, until all stubs are
connected.

The probability mass function and the distribution function, of the
nodal degree, are denoted by
    \begin{equation}
    \label{kansen}
    \prob(D_1=j)=f_j,\quad j=0,1,2,\ldots, \quad \mbox{and} \quad
    F(x)=\sum_{j=0}^{\lfloor x \rfloor} f_j,
    \end{equation}
where $\lfloor x \rfloor$ is the largest integer smaller than or
equal to $x$. Our main assumption will be that
    \begin{equation}
    \label{distribution}
    1-F(x)=x^{-\tau+1}L(x),
    \end{equation}
where $\tau \in (1,2)$ and $L$ is slowly varying at infinity. This means
that the random variables $D_j$ obey a power law with infinite mean. The factor
$L$ is meant to generalize the model.

\subsection{Main results}
\label{sec-res}

We define the graph distance $H_{\sN}$ between
the nodes $1$ and $2$ as the minimum number of edges that form a
path from $1$ to $2$, where, by convention, this distance
equals $\infty$ if $1$ and $2$ are not connected.
Observe that the distance between two randomly chosen
nodes is equal in distribution to $H_{\sN}$,
because the nodes are exchangeable.

We will present in this paper two separate theorems
for the case $\tau\in (1,2)$ .
In the first theorem we take the sequence $D_1,D_2,\ldots,D_{\sN}$ an i.i.d.\,
sequence with distribution $F$, satisfying
(\ref{distribution}), with $\tau \in (1,2)$. The result
is that the graph distance or hopcount converges in distribution to
a limit random variable with mass $p=p_\tau$, $1-p,$ on the values
$2$, $3$, respectively.

One might argue that including degrees larger than $N-1$ is
artificial in a network with $N$ nodes. In fact, in many real networks,
the degree is bounded by a physical constant. Therefore we also consider
the case where the degrees are conditioned to be smaller than $N^{\alpha}$, where
$\alpha$ is an arbitrary positive number. Of course, we cannot condition
on the degrees to be at most $M$, where $M$ is fixed and independent
on $N$, since in this case, the degrees are uniformly bounded, and this case is
treated in \cite{HHV03}. Therefore, we consider cases where the degrees are
conditioned to be at most a given power of $N$.

The result with conditioned degrees appears in the second theorem below.
It turns out that for $\alpha >1/(\tau-1)$, the conditioning has no influence
in the sense that the limit random variable is the same as that for the
unconditioned case. This is not so strange, since the maximal degree is
of order $N^{1/(\tau-1)}$, so that the conditioning does nothing
in this case. However, for $1/\tau< \alpha <1/(\tau-1)$, the graph
distance converges to a degenerate limit random variable with mass
$1$ on the value $3$. It would be of interest to extend the possible conditioning
schemes, but we will not elaborate further on it here.

\begin{theorem}
\label{thm-tau(1,2)-uncond} Fix $\tau \in (1,2)$ in
(\ref{distribution}) and let $D_1,D_2,\ldots,D_{\sN}$ denote an i.i.d.\
sequence with common distribution $F$. Then,
\begin{equation}
\label{limit law}
\lim_{N\to \infty} \prob(H_{\sN}=2)
=1-\lim_{N\to \infty} \prob(H_{\sN} =3)=p,
\end{equation}
where $p=p_{\tau}\in (0,1)$.
\end{theorem}

In the theorem below we will write $D^{\smallsup{N}}$ for the random variable
$D$ conditioned on $D<N^{\alpha}$. Thus,
    \eq
    \label{truncate}
    \prob(D^{\smallsup{N}}=k)=\frac{f_k}{\prob(D<N^\alpha)}, \qquad 0\leq k< N^{\alpha}.
    \en

\begin{theorem}
\label{thm-tau(1,2)-cond} Fix $\tau \in (1,2)$ in
(\ref{distribution}) and let
$D_1^{\smallsup{N}},D_2^{\smallsup{N}},\ldots, D_{\sN}^{\smallsup{N}}$
be a sequence of i.i.d.\  copies of $D^{\smallsup{N}}$.
\begin{enumerate}
\item[{\rm (i)}] If $1/\tau<\alpha<1/(\tau-1)$, then
    \begin{equation}
    \label{limit law2i}
    \lim_{N\to \infty}
    \prob(H_{\sN} =3)=1.
\end{equation}
\item[{\rm (ii)}] If $\alpha>1/(\tau-1)$, then
    \begin{equation}
    \label{limit law2ii}
    \lim_{N\to \infty}
    \prob(H_{\sN}=2)=1-\lim_{N\to \infty}
    \prob(H_{\sN}=3)=p,
    \end{equation}
where $p=p_{\tau}\in (0,1)$.
\end{enumerate}
\end{theorem}

\noindent
{\bf Remark:} For $\alpha<1/\tau$, we have reasons to believe that
$\prob(H_{\sN}>3)$ remains uniformly positive when $N\rightarrow \infty$.
A heuristic for this fact is given in the next section.

\subsection{Heuristics}
\label{sec-heuristics}
When $\tau\in (1,2)$, we consider two different cases.
In Theorem \ref{thm-tau(1,2)-uncond}, the degrees are
not conditioned, while in Theorem \ref{thm-tau(1,2)-cond}
we condition on each node having a degree smaller than
$N^{\alpha}$. We will now give a heuristic explanation of
our results.

In two previous papers \cite{HHV03, HHZ04b}, we have
treated the cases $\tau \in (2,3)$ and $\tau >3$. In these cases, it follows
that the probability mass function $\{f_j\}$ introduced in (\ref{kansen})
has a finite mean, and the the number of nodes on graph
distance $k$ from node $1$ can be coupled to the number of individuals in the
$k^{\emph{th}}$ generation of a branching process with offspring distribution
$\{g_j\}$ given by
    \eq
    g_j = \frac{j+1}{\mu} f_j,
    \en
where $\mu=\expec[D_1]$. For $\tau \in (1,2)$, as we are currently investigating,
we have $\mu=\infty$, and the branching process used in \cite{HHV03, HHZ04b}
does not exist.

When we do not condition on $D_j$ being smaller than
$N^{\alpha}$, then $L_{\sN}$ is the i.i.d.\
sum of $N$ random variables $D_1,D_2,\ldots,D_{\sN}$, with infinite
mean. It is well known that in this case the bulk of the
contribution to $L_{\sN}\sim N^{1/(\tau-1)}$ comes from a {\it finite} number of nodes
which have giant degrees (the so-called {\it giant nodes}). A basic fact
in the configuration model is that two sets of stubs of sizes
$n$ and $m$ are connected with high probability when $nm$ is at least of order
$L_{\sN}$. Since the giant nodes have degree roughly $N^{1/(\tau-1)}$, which is
much larger than $\sqrt{L_{\sN}}$, they are all attached to each other,
thus forming a complete graph of giant nodes. Each stub is with
probability close to 1 attached to a giant node, and therefore,
the distance between any two nodes is, with large probability, at
most 3. In fact this distance equals 2 when the two nodes are
connected to the {\it same} giant node, and is 3 otherwise.

When we truncate the distribution as in (\ref{truncate}),
with $\alpha>1/(\tau-1)$, we  hardly change anything
since without truncation with probability $1-o(1)$ all
degrees are below $N^\alpha$.
On the other hand, if $\alpha<1/(\tau-1)$ when, with truncation,
the largest nodes have degree of order $N^{\alpha}$,
and $L_{\sN}\sim N^{1+\alpha(2-\tau)}$. Again, the bulk of
the total degree $L_{\sN}$ comes from nodes with degree
of the order $N^{\alpha}$, so that now these are the giant nodes.
Hence, for $1/\tau<\alpha < 1/(\tau-1)$, the largest nodes
have degrees much larger than $\sqrt{L_{\sN}}$,
and thus, with probability $1-o(1)$, still constitute a
complete graph. The number of giant nodes converges to infinity,
as $N\to \infty$. Therefore, the probability that two arbitrary
nodes are connected to the {\it same} giant node converges to 0.
Therefore, the hopcount equals 3 with probability converging to 1.
If $\alpha<1/\tau$, then the giant nodes no longer constitute a
complete graph suggesting that the resulting hopcount
can be greater than $3$. It remains unclear to us whether the hopcount
converges to a {\it single} value (as in Theorem \ref{thm-tau(1,2)-cond}),
or to more than one possible value (as in Theorem \ref{thm-tau(1,2)-uncond}).
We do expect that the hopcount remains uniformly bounded.

The proof in this paper is based on detailed asymptotics of the sum
of $N$ i.i.d.\ random variables with infinite mean, as well as on the
scaling of the order statistics of such random variables. The scaling
of these order statistics is crucial in the definition of the giant nodes
which are described above. The above considerations are the basic idea in
the proof of Theorem \ref{thm-tau(1,2)-uncond}. In the proof of Theorem
\ref{thm-tau(1,2)-cond}, we need to investigate what the conditioning does
to the scaling of both the total degree $L_{\sN}$, as well as to the largest
degrees.

\subsection{Related work}
\label{sec-related}
The above model is a variation of the configuration model. In the
usual configuration model one often starts from a given {\it deterministic} degree sequence.
In our model, the degree sequence is an i.i.d. sequence
$D_1, \ldots, D_{\sN}$ with distribution equal to a power law. The
reason for this choice is that we are interested in models for
which all nodes are exchangeable, and this is not the case when
the degrees are fixed. The study of this variation of the configuration model was started
in \cite{NSW00} for the case $\tau>3$ and studied by Norros and Reittu \cite{norros} in case
$\tau\in (2,3)$.

For a complete survey to complex networks, power law degree sequences and
random graph models for such networks, see \cite{AB02} and \cite{Newm03}. There
a heuristic is given why the hopcount scales proportionally to $\log N$, which is
originally from \cite{NSW00}.
The argument uses a variation of the power law degree model, namely, a model where
an exponential cut off is present. An example of such a degree distribution is
    \eq
    f_k=C k^{-\tau} e^{-k/\kappa}
    \en
for some $\kappa>0$. The size of $\kappa$ indicates up to what
degree the power law still holds, and where the exponential
cut off starts to set in. The above model is treated in \cite{HHV03}
for any $\kappa<\infty$, but, for $\kappa=\infty$, falls within
the regimes where $\tau\in (2,3)$ in \cite{HHZ04b} and within the
regime in this paper for $\tau\in (1,2)$. In \cite{NSW00}, the authors
conclude that since the limit as $\kappa\rightarrow \infty$ does not
seem to converge, the `average distance is not well-defined when $\kappa<3$'.
In this paper, as well as in \cite{HHZ04b}, we show that the average distance
{\it is} well-defined, but it scales differently from the case where
$\tau>3$.

In the paper \cite{HHZ04c}, we give a survey to the results for the hopcount
in the three different regimes $\tau\in (1,2)$, $\tau\in (2,3)$ and $\tau>3$.
There, we also prove results for the connectivity properties of the random graph
in these cases. These results assume that the expected degree is larger than 2.
This is always the case when $\tau\in(1,2)$, and stronger results have been shown
there. We prove that the largest connected component has size $N(1+o(1))$ with
probability converging to one.
When $\tau\in (1,\frac 32)$, we can even prove that with large probability, the graph
is connected. When $\tau>\frac 32$, this is not true, and we investigate the structure
of the remaining `dust' that does not belong to the largest connected component.
In the analysis, we will make use of the results obtained in this paper for $\tau\in(1,2)$.
For instance, it will be crucial that the probability that two arbitrary nodes are connected
converges to 1.

There is substantial related work on the configuration model
for the case $\tau \in (2,3)$ and $\tau>3$. References are included in the paper
\cite{HHZ04b} for  the case $\tau \in (2,3)$ and in \cite{HHV03} for  $\tau>3$.
We again refer to the references in \cite{HHZ04c} and \cite{AB02, Newm03}
for more details. The graph distance for $\tau \in (1,2)$, that we study here,
has to our best knowledge not been studied before. Values of $\tau \in (1,2)$
have been observed in networks of e-mail messages and networks
where the nodes consist of software packages (see \cite[Table II]{Newm03}),
for which our configuration model with $\tau\in (1,2)$ can possibly give
a good model.

In \cite{ACL01}, random graphs are considered with a degree
sequence that is {\it precisely} equal to a power law, meaning
that the number of nodes with degree
$k$ is precisely proportional to $k^{-\tau}$. Aiello {\em et al.} \cite{ACL01}
show that the largest connected component is of the order of the size of the graph
when $\tau<\tau_0=3.47875\ldots$, where $\tau_0$ is the solution of
$\zeta(\tau-2)-2\zeta(\tau-1)=0$, and where $\zeta$ is the Riemann Zeta
function. When $\tau>\tau_0$, the largest connected component is
of smaller order than the size of the graph and more precise
bounds are given for the largest connected component. When
$\tau\in (1,2)$, the graph is with high probability connected. The proofs
of these facts use couplings with branching processes and strengthen
previous results due to Molloy and Reed \cite{MR95,MR98}.
See also \cite{ACL01} for a history of the problem and references predating
\cite{MR95,MR98}. The problem of distances in the configuration model with
$\tau \in (1,2)$ has, up to our best knowledge, not
been addressed. See \cite{ACL01c} for an introduction to the
mathematical results of various models for complex networks
(also called massive graphs), as well as a detailed account of the results
in \cite{ACL01}.

A detailed account for a related model can be found in
\cite{CL02a} and \cite{CL02b}, where links
between nodes $i$ and $j$ are present with probability
equal to $w_iw_j/\sum_l w_l$ for some `expected degree vector'
$w=(w_1, \ldots, w_N)$. Chung and Lu \cite{CL02a} show that when $w_i$ is
proportional to $i^{-{\frac 1{\tau-1}}}$ the average distance
between pairs of nodes is $\log_{\nu}N(1+o(1))$ when $\tau>3$, and
$2\frac{\log\log N}{|\log(\tau-2)|}(1+o(1))$ when $\tau\in (2,3)$.
In their model, also $\tau<1$ is possible, and in this case,
similarly to $\tau\in (1,\frac 32)$ in our paper, the graph is
connected with high probability.

The difference between this model and ours is that the nodes are
not exchangeable in \cite{CL02a}, but the observed phenomena
are similar. This result can be understood as follows. Firstly,
the actual degree vector in \cite{CL02a} should be close to the
expected degree vector. Secondly, for the expected degree vector, we can compute that
the number of components for which the degree is less than or equal to $k$ equals
$$
|\{i: w_i\leq k\}|\propto |\{i: i^{-\frac{1}{\tau-1}}\leq k\}|\approx k^{-\tau+1}.
$$
Thus, one expects that the number of nodes with degree at most $k$ decreases as
$k^{-\tau+1}$, similarly as in our model. In \cite{CL02b}, Chung and Lu study
the sizes of the connected components in the above model. The advantage of working with
the `expected degree model' is that different links are present independently
of each other, with makes this model closer to the random graph $G(p,N)$.

\subsection{Organization of the paper}
\label{sec-organization}
The main body of the paper consists of the proofs of Theorem \ref{thm-tau(1,2)-uncond}
in Section \ref{sec-pftau(1,2)-uncond} and the proof
of Theorem \ref{thm-tau(1,2)-cond} in Section \ref{sec-pftau(1,2)-cond}. Both proofs
contain a technical lemma and in order to make the argument more transparent, we
have postponed the proofs of these lemmas to the appendix.
Section \ref{sec-sim} contains simulation results and some conclusions.


\section{Proof of Theorem~\ref{thm-tau(1,2)-uncond}}
\label{sec-pftau(1,2)-uncond}

In this section, we prove Theorem~\ref{thm-tau(1,2)-uncond},
which states that the hopcount between two arbitrary nodes has, with
probability $1-o(1)$, a non-trivial distribution on $2$ and $3$.

We will use an auxiliary lemma, which is a modification of the
extreme value theorem for the $k$ largest degrees, $k\in\NN$.
We introduce
    $$
    D_{\smallsup{1}}\le D_{\smallsup{2}}\le\dots\le D_{\smallsup{N}},
    $$
to be the order statistics of $D_1,\dots,D_{\sN}$, so that
$D_{{\smallsup{1}}}=\min \{D_1,\ldots,D_{\sN}\}$, $D_{\smallsup{2}}$
is the second smallest degree, etc.
Let $(u_{\sN})$ be an increasing sequence such that
    \begin{equation}
    \label{dnz11}
    \lim\limits_{N\to\infty}N\left(1-F(u_{\sN})\right)=1.
    \end{equation}
It is well known that the order statistics of the degrees, as well as the total degree,
are governed by $u_{\sN}$ in the case that $\tau\in (1,2)$. The following lemma
shows this in some detail.
\begin{lemma}
\label{Ldnz1}
\par
$(a)$ for any $k\in\NN$,
    \begin{equation}
    \label{dnz67_3}
    \left(\frac{D_{\smallsup{N}}}{u_{\sN}},\dots,
    \frac{D_{\smallsup{N-k+1}}}{u_{\sN}}\right)
    \longrightarrow\left(\xi_1,\dots,\xi_k\right), \mbox{ in
    distribution, as }N\to\infty,
    \end{equation}
where $(\xi_1,\dots,\xi_k)$ is a random
vector with marginals in $(0,\infty)$ and with joint distribution
function given, for any tuple $0<y_k<\dots<y_1<\infty$, by
    \begin{eqnarray}
    \label{dnz67_5}
    &&\prob\left(\xi_1<y_1,\dots,\xi_k<y_k\right)\\
    &&\qquad=\sumu\limits_{0\le r_1\le\dots\le r_k<k}
    \frac{y_1^{(1-\tau)r_1}}{r_1!}
    \frac{(y_2^{1-\tau}-y_1^{1-\tau})^{r_2}}{r_2!} \dots
    \frac{(y_k^{1-\tau}-y_{k-1}^{1-\tau})^{r_k}}{r_k!}
    e^{-y_k^{1-\tau}}.\nn
    \end{eqnarray}
Moreover,
\begin{equation}
\label{dnz67_6}
\xi_k\to0,\mbox{ in probability, as }k\to\infty.
\end{equation}
\noindent
$(b)$
$$
\frac{L_{\sN}}{u_{\sN}}\longrightarrow\eta,
\mbox{ in distribution, as }N\to\infty,$$
where $\eta$ is a random variable on $(0,\infty)$.
\end{lemma}

\begin{proof} For part $(a)$, we take $\rho_i=y_i^{1-\tau}$ and
$u_{\sN}(\rho_i)=y_iu_{\sN}$, $i\in \{1,\dots,k\}$. Since
$u_{\sN}=L^{'}(N)N^{\frac{1}{\tau-1}}$ (see e.g., \cite{EKM97})
for some slowly varying function $L^{'}(N)$, it follows
from (\ref{dnz11}) that,
    $$
    \lim\limits_{N\to\infty}N\left(1-F(u_{\sN}(\rho_i))\right)=\rho_i,
    \qquad i\in\{1,\dots,k\}.
    $$
Hence by \cite[Theorem 4.2.6 and (4.2.4)]{EKM97}, we have
    \begin{eqnarray*}
    &&\lim_{N\to\infty}
    \prob\left(D_{\smallsup{N}}<u_{\sN}(\rho_1),\dots,D_{\smallsup{N-k+1}}<u_{\sN}(\rho_k)\right)\\
    &&\qquad=\sumu\limits_{0\le r_1\le\dots\le r_k<k}
    \frac{{\rho_1}^{r_1}}{r_1!} \frac{(\rho_2-\rho_1)^{r_2}}{r_2!}
    \dots \frac{(\rho_k-\rho_{k-1})^{r_k}}{r_k!}
    e^{-\rho_k}\\
    &&\qquad=\sumu\limits_{0\le r_1\le\dots\le r_k<k}
    \frac{y_1^{(1-\tau)r_1}}{r_1!}
    \frac{(y_2^{1-\tau}-y_1^{1-\tau})^{r_2}}{r_2!} \dots
    \frac{(y_k^{1-\tau}-y_{k-1}^{1-\tau})^{r_k}}{r_k!}
    e^{-y_k^{1-\tau}}.
    \end{eqnarray*}
We compute now the marginal distribution of $\xi_k$, for $k\ge1$.
For any $x>0$, due to~(\ref{dnz67_5}), we have
$$
\begin{array}{rl}
\prob(\xi_k<x)
&=\prob(\xi_1<\infty,\dots,\xi_{k-1}<\infty,\xi_k<x)\\
&=\lim\limits_{y_1,\dots,y_{k-1}\to\infty}
\sumu\limits_{0\le r_1\le\dots\le r_{k-1}\le r_k<k}
e^{-x^{1-\tau}}\prod_{i=1}^k
\frac{\big(y_i^{(1-\tau)r_i}-y_{i-1}^{(1-\tau)r_{i-1}}\big)}{r_i!}
\\
&=\sumu\limits_{r=0}^{k-1}
\frac{x^{(1-\tau) r}}{r!}e^{-x^{1-\tau}}\to 1,\mbox{ as }k\to\infty,
\end{array}
$$
where in the middle expression, we write $y_0=0$ and $y_k=x$.
Hence we have~(\ref{dnz67_6}).

Part $(b)$ follows since $L_{\sN}=D_1+\dots+D_{\sN}$ is in the domain of attraction of a stable law
(\cite[Corollary 2, XVII.5, p.~578]{fellerb}).
\end{proof}
\qed

\medskip

We need some additional notation.
In this section (Section \ref{sec-pftau(1,2)-uncond}) we define the {\em giant
nodes} as the $k_\varepsilon$ largest nodes, i.e., those nodes with degrees
$D_{\smallsup{N}},\dots,D_{\smallsup{N-k_\varepsilon+1}},$ where
$k_\varepsilon$ is some function of $\varepsilon$, to
be chosen below. We define
\begin{equation}
\label{def-AeN}
\AeN=\BeN\cap\CeN\cap\DeN,
\end{equation}
where
\begin{itemize}
\item[(i)]
$\BeN$ is the event that the stubs of node 1 and node 2 are
attached exclusively to stubs of giant nodes.
\item[(ii)]$\CeN$ is the event that any
two giant nodes are attached to each other; and
\item[(iii)] $\DeN$ is defined as
$$
\DeN=\left\{D_1\le b_{\sD,\eps}\cap D_2\le b_{\sD,\eps}\right\},
$$
where  $b_{\sD,\eps} =\min\{k: 1-F(k)<\eps/8\}$, so that $b_{\sD,\eps}= \eps^{-1/(\tau-1)(1+o(1))}$.
\end{itemize}
The sets $\BeN$ and $\CeN$ depend on the integer $k_{\eps}$, which we will take
to be large for $\eps$ small, and will be defined now.
The choice of the index $k_{\eps}$ is rather technical, and depends on
the distributional limits defined in Lemma \ref{Ldnz1}.
Since $L_{\sN}/u_{\sN}$ converges in distribution to the random variable $\eta$, with support
$(0,\infty)$, we can find $a_{\eta, \eps}$, such that
    \begin{equation}
    \label{defa3N}
    \prob(L_{\sN}< a_{\eta, \eps} u_{\sN})< \eps/36, \quad \forall N.
    \end{equation}
This follows since convergence in distribution implies tightness of the sequence $L_{\sN}/u_{\sN}$ (\cite[p.~9]{bill68}),
so that we can find a closed subinterval $[a,b]\subset (0,\infty)$, with
    $$
    \prob(L_{\sN}/u_{\sN}\in [a,b])>1-\eps,\quad \forall N.
    $$
The definition of $b_{\xi_{k_\eps}}$ is rather involved. It depends on $\eps$, the quantile $b_{\sD,\eps}$, the value $a_{\eta, \eps}$
defined above and the value of $\tau\in (1,2)$ and reads
    \begin{equation}
    \label{defa4N}
    b_{\xi_{k_\eps}}
    =\left(\frac{\eps^2a_{\eta, \eps}}{2304b_{\sD,\eps}}\right)
    ^{\frac{1}{2-\tau}},
    \end{equation}
where the peculiar integer $2304$ is the product of $8^2$ and $36$.
Given $b_{\xi_{k_\eps}}$, we take $k_\eps$ equal to the minimal value such that
    \begin{equation}
    \label{defkeps}
    \prob(\xi_{k_\eps}\geq b_{\xi_{k_\eps}}/2)\le\varepsilon/72.
    \end{equation}
This we can do due to~(\ref{dnz67_6}). We have now defined the constants
that we will use in the proof, and we will next show that the probability of
$\AeN^c$ is at most $\eps$:
    \begin{lemma}
    \label{technlemma1}
    For each $\eps>0$, there exists $N_\eps$, such that
    \begin{equation}
    \label{afschattingAeN}
    \prob(\AeN^c)< \eps, \quad N\geq N_\eps.
    \end{equation}
    \end{lemma}
The proof of this Lemma is rather technical and can be found in the appendix.
We will now complete the proof of Theorem \ref{thm-tau(1,2)-uncond} subject to
Lemma \ref{technlemma1}.
\medskip

\noindent{\bf Proof of Theorem \ref{thm-tau(1,2)-uncond}.}
The proof consist of two parts. The event $\AeN$,
implies the event $\{H_{\sN}\le3\}$, so that $\prob(\AeN^c)<\varepsilon$
induces that $\{H_{\sN}\le3\}$ with probability at least $1-\eps$.

In the first part we show that $\prob\left(\{H_{\sN}=1\}\cap \AeN \right)<\eps$.
In the second part we prove that
$$
\lim_{N\to\infty}\prob\left(H_{\sN}=2\right)
=\lim\limits_{\varepsilon\to 0}\lim_{N\to\infty}
\prob\left(\{H_{\sN}=2\}\cap\DeN\,|\,\BeN\right)=p,
$$
for the events $\BeN,\DeN$ defined above and some $0<p<1$.
Since $\eps$ is arbitrary positive, the above statements yield the content of the theorem.

We first prove that $\prob\left(\{H_{\sN}=1\}\cap \AeN \right)<\eps$
for sufficiently large $N$.
The event $\{H_{\sN}=1\}$ occurs {\it iff} at
least one stub of node $1$ connects to a stub of node $2$. For
$j\le D_1$,  we denote by $\{[1.j]\to[2]\}$ the event
that $j$-th stub of node $1$ connects to a stub of node $2$. Then,
with $\prob_{\sN}$ the conditional probability given the degrees $D_1,D_2,\ldots,D_N$,
\begin{eqnarray}
\label{hop=1}
\prob(H_{\sN}=1,\,A_{\varepsilon, \sN}) &\le&
\expec\left[
\sum_{j=1}^{D_1}
\prob_{\sN}(\{[1.j]\to[2]\},\,A_{\varepsilon, \sN})\right]\nonumber\\
&\le&\expec\left[ \sumu\limits_{j=1}^{D_1}
\frac{D_2}{L_{\sN}-1}\indic\{A_{\varepsilon, \sN}\}\right]
\le\frac{b_{\sD,\eps}^2}{N-1} <\varepsilon,
\end{eqnarray}
for large enough $N$, since $L_{\sN}\ge N$.

We prove now that
$\lim\limits_{N\to\infty}\prob\left(H_{\sN}=2\right)=p$, for some
$0<p<1$.
Since by definition for any $\eps>0$,
$$
\max\{\prob(\BeN^c),\prob(\DeN^c)\}\le \prob(\AeN^c)\le\eps,
$$
we have that
\begin{eqnarray*}
&&|\prob(H_{\sN}=2)-
\prob\left(\{H_{\sN}=2\}\cap\DeN\,|\,\BeN \right)|\\
&&\quad\leq
\left|\prob(H_{\sN}=2)\left(1-\frac{1}{\prob(\BeN)}\right)\right|
+\left|\frac{\prob(H_{\sN}=2)-\prob(\{H_{\sN}=2\}\cap\DeN\cap \BeN)}{\prob(\BeN)}\right|\\
&&\quad\leq
\frac{2{\prob(\BeN^c)+\prob(\DeN^c)}}{\prob(\BeN)}\le \frac{3\eps}{1-\eps},
\end{eqnarray*}
uniformly in $N$, for $N$ sufficiently large.
If we show that
$$
\lim\limits_{N\to\infty} \prob\left(\{H_{\sN}=2\}\cap\DeN\,|\,\BeN
\right)=p_{\tau,\varepsilon},
$$
then there exists a double limit
$$
p_\tau=\lim\limits_{\varepsilon\to0}\lim\limits_{N\to\infty}
\prob\left(\{H_{\sN}=2\}\cap\DeN\,|\,\BeN \right)
=\lim\limits_{N\to\infty}\prob\left(H_{\sN}=2\right).
$$
Moreover, if
we can bound $p_{\tau,\varepsilon}$ from $0$ and $1$, uniformly in $\varepsilon$,
for $\varepsilon$ small enough, then we also obtain that $0<p_\tau<1$.

First, to prove the existence of
    \begin{equation}
    \label{dnz67_7}
    \lim\limits_{N\to\infty} \prob\left(\{H_{\sN}=2\}\cap\DeN\,|\,\BeN\right)
    =\lim\limits_{N\to\infty}\expec\left(\prob_{\sN}\left(\{H_{\sN}=2\}\cap\DeN\,|\,\BeN\right)\right)
    =p_{\tau,\varepsilon},
    \end{equation}
we will show that
$\prob_{\sN}\left(\{H_{\sN}=2\}\cap\DeN\,|\,\BeN \right)$
is a continuous function of the vector
    $$
    \bar D_{k_\eps}=\left(\frac{D_{\smallsup{N}}}{u_{\sN}},\dots
    \frac{D_{\smallsup{N-k_\varepsilon+1}}}{u_{\sN}},
    \frac{1}{u_{\sN}}\right).
    $$
This vector, due to~(\ref{dnz67_3}),
converges in distribution to $\left(\xi_1,\dots,\xi_{k_\varepsilon},0\right)$.
Hence, by the continuous mapping theorem
\cite[Theorem 5.1, p.~30]{bill68}, we have
the existence of the limit~(\ref{dnz67_7}).
We now prove the claimed continuity.

We will show an even stronger statement, namely, that
$\prob_{\sN}\left(\{H_{\sN}=2\}\cap\DeN\,|\,\BeN \right)$ is the
ratio of two polynomials of the components of the vector $\bar D_{k_\eps}$,
where the polynomial in the denominator is strictly positive.
Indeed, the hopcount between nodes $1$ and $2$ is $2$
iff both nodes are connected to the same giant node.
For any $0\le i\le D_1$, $0\le j\le D_2$ and $0\le k<k_\eps$
let ${\cal A}_{i,j,k}$ be the event that both the $i^{\rm th}$ stub of node 1
and the $j^{\rm th}$ stub of node 2 are connected to the node with the $(N-k)^{\rm th}$
largest degree.
Then,
    $$
    \prob_{\sN}\left(\{H_{\sN}=2\}\cap\DeN\,|\,\BeN \right)
    =\prob_{\sN}\left(\bigcup\limits_{i=1}^{D_1}
    \bigcup\limits_{j=1}^{D_2}
    \bigcup\limits_{k=0}^{k_\eps-1}
    {\cal A}_{i,j,k}\,\Big|\,\BeN\right),
    $$
where the r.h.s.\ can be written by inclusion-exclusion
formula, as a linear combination of terms
    \eq
    \label{inters}
    \prob_{\sN}\left({\cal A}_{i_1,j_1,k_1}\cap\cdots\cap{\cal A}_{i_n,j_n,k_n}\,|\,\BeN\right).
    \en
It is not difficult to see that these probabilities
are ratios of polynomials. For example,
    \eqalign
    \label{dnz67_8}
    \prob_{\sN}\left({\cal A}_{i,j,k}\,|\,\BeN\right)
    &=\frac{D_{\smallsup{N-k}}(D_{\smallsup{N-k}}-1)}{(D_{\smallsup{N-k_{\eps}+1}}+\dots+D_{\smallsup{N}})
    (D_{\smallsup{N-k_{\eps}+1}}+\dots+D_{\smallsup{N}}-1)}\\
    &=\frac{\frac{D_{\smallsup{N-k}}}{u_{\sN}}(\frac{D_{\smallsup{N-k}}}{u_{\sN}}-\frac{1}{u_{\sN}})}
    {(\frac{D_{\smallsup{N-k_{\eps}+1}}}{u_{\sN}}+\dots+\frac{D_{\smallsup{N}}}{u_{\sN}})
    (\frac{D_{\smallsup{N-k_{\eps}+1}}}{u_{\sN}}+\dots+\frac{D_{\smallsup{N}}}{u_{\sN}}-\frac{1}{u_{\sN}})}.
    \nn
    \enalign
Similar arguments hold for general terms of the form in (\ref{inters}).
Hence, $\prob_{\sN}\left(\{H_{\sN}=2\}\cap\DeN\,|\,\BeN \right)$ itself
can be written as a ratio of two polynomials where the polynomial in the
denominator is strictly positive. Therefore, the limit~(\ref{dnz67_7}) exists.

We finally bound $p_\eps$ from $0$ and $1$
uniformly in $\eps$, for any $\eps<1/2$.
Since the hopcount between nodes $1$ and $2$ is $2$,
given $\BeN$,
if they are both connected to the node with largest degree,
    $$
    \prob(\{H_{\sN}=2\}\cap\DeN\,|\,\BeN)
    \ge
    \expec\left[\prob_{\sN}\left({\cal A}_{\sss 1,1,(N)}\,|\,\BeN\right)\right],
    $$
and by~(\ref{dnz67_8}) we have
    $$
    \begin{array}{rl}
    p_{\tau,\eps}&=\lim\limits_{N\to\infty}
    \prob\left(\{H_{\sN}=2\}\cap\DeN\,|\,\BeN \right)
    \ge\lim\limits_{N\to\infty}
    \expec\left[\frac{D_{\smallsup{N}}(D_{\smallsup{N}}-1)}{(D_{\smallsup{N}}+\cdots+D_{\smallsup{N-k_\varepsilon+1}}-1)^2}\right]\\
    &=\expec\left[\left(\frac{\xi_1}{\xi_1+\cdots+\xi_{k\varepsilon}}\right)^2\right]
    \ge\expec\left[\left(\frac{\xi_1}{\eta}\right)^2\right].
    \end{array}
    $$
On the other hand,
the hopcount between nodes $1$ and $2$ is at most $3$, given $\BeN$
when all stubs of the node $1$ are connected to the node with largest degree,
and all stubs of the node $2$ are connected to the node with the one but
largest degree. Hence, for any $\varepsilon<1/2$ and similarly to~(\ref{dnz67_8}), we have
    $$
    \begin{array}{rl}
    p_{\tau,\varepsilon}&=\lim\limits_{N\to\infty}
    \prob\left(\{H_{\sN}=2\}\cap\DeN\,|\,\BeN \right)
    \le1-\lim\limits_{N\to\infty}
    \prob\left(\{H_{\sN}>2\}\cap\DeN\,|\,\BeN \right) \\
    &\le1-\lim\limits_{N\to\infty}
    \prob\left(\{H_{\sN}>2\}\cap
    D_{{\sss \frac 1 2} \kern-.05em,\kern .1em \sN}\,|\,\BeN \right)\\
    &\le1-\lim\limits_{N\to\infty}
    \expec\left[\left(\prod\limits_{i=0}^{D_1}
    \frac{D_{(N)}-2i}{D_{\smallsup{N}}+\cdots+D_{(N-k_\varepsilon+1)}-D_1}
    \right)
    \left(\prod\limits_{i=0}^{D_2}
    \frac{D_{(N-1)}-2i}{D_{\smallsup{N}}+\cdots+D_{(N-k_\varepsilon+1)}-D_2}
    \right)
    \indic\{D_{\frac 1 2 \kern-.05em,\kern .1em \sN}\}
    \right]\\
    &\le1-\frac{1}{2}\expec\left[\left(\frac{\xi_1}{\xi_1+\cdots+\xi_{k_\varepsilon}}\right)^{b_{\sD,1/2}}
    \left(\frac{\xi_2}{\xi_1+\cdots+\xi_{k_\varepsilon}}\right)^{b_{\sD,1/2}}\right]
    \le1-\frac{1}{2}\expec\left[\left(\frac{\xi_1\xi_2}{\xi^2}\right)^{b_{\sD,1/2}}\right].
    \end{array}
    $$
Both the upper and lower bound are strictly positive, independently of $\eps.$
Hence, for any $\varepsilon<1/2$,
the quantity $p_{\tau,\eps}$ is bounded away from $0$ and $1$, where the bounds
are {\it independent of $\eps$}, and thus $0<p=p_{\tau}<1$.
\qed

This completes the proof of Theorem \ref{thm-tau(1,2)-uncond} subject to
Lemma \ref{technlemma1}.


\section{Proof of Theorem \ref{thm-tau(1,2)-cond}}
\label{sec-pftau(1,2)-cond}

In Theorem \ref{thm-tau(1,2)-cond}, we consider the hopcount in the configuration model
with degrees an i.i.d. sequence with distribution given by (\ref{truncate}),
where $D$ has distribution $F$ satisfying (\ref{distribution}).
We distinguish two cases: (i) $1/\tau<\alpha <1/(\tau-1)$ and (ii) $\alpha > 1/(\tau-1)$.

We first prove part (ii), which states that the limit distribution of
$H_{\sN}$ is a mixed distribution with probability mass
$p$ on  $2$ and probability mass $1-p$ on $3$, for some $0<p<1$.
Part (ii) of Theorem \ref{thm-tau(1,2)-cond} is almost immediate
from Theorem \ref{thm-tau(1,2)-uncond}.
As before we denote by $D_1,D_2,\ldots,D_{\sN}$ the i.i.d. sequence without conditioning, then $\prob(\cup_{i=1}^N
\{D_i>N^{\alpha}\})$, which is the probability that for at least one index $i\in \{1,2,\ldots,N\}$ the degree $D_i$
exceeds $N^{\alpha}$, is bounded by
$$
\sum_{i=1}^N \prob(D_i>N^{\alpha})= N\prob(D_1>N^{\alpha})=N^{1+\alpha(1-\tau)}L(N
)=N^{-\eps},
$$
for some positive $\eps$, because
$\alpha > 1/(\tau-1)$. We can therefore couple the i.i.d.\ sequence
$\vec D^{\smallsup{N}}=(D_1^{\smallsup{N}},D_2^{\smallsup{N}},\ldots,
D_{\sN}^{\smallsup{N}})$ to the sequence $\vec D=(D_1,D_2,\ldots,D_{\sN})$,
where the probability of a miscoupling, i.e., a coupling such that
$\vec D^{\smallsup{N}}\neq \vec D,$ is at most $N^{-\eps}.$
Therefore, the result of Theorem \ref{thm-tau(1,2)-uncond} carries over to case (ii)
in Theorem \ref{thm-tau(1,2)-cond}.

\medskip
We now turn to case (i) in Theorem \ref{thm-tau(1,2)-cond}.
We must prove that if we condition the degrees to be smaller than
$N^{\alpha}$ with $1/\tau<\alpha<1/(\tau-1)$, then the graph distance between
two arbitrary nodes is $3$ with probability equal to $1-o(1)$.
We define
    $$
    L^{\smallsup{N}}_{\sN}=\sumu\limits_{n=1}^ND_n^{\smallsup{N}},
    $$
to be the total degree of the conditioned model.
The steps in the proof Theorem \ref{thm-tau(1,2)-cond}(i) are
identical to those in the proof of Theorem \ref{thm-tau(1,2)-uncond},
however, the details differ considerably.
We take an arbitrary $\varepsilon>0$ and define an event
$\GeN$ such that
$\prob(\GeN^c)<\varepsilon$, and prove subsequently that
for large enough $N$,
$\prob(\{H_{\sN}=1\}\cap\,\GeN)<\varepsilon$,
$\prob(\{H_{\sN}=2\}\cap\,\GeN)<\varepsilon$, and that
$\prob(\{H_{\sN}\le3\}\cap\,\GeN)\ge1-\varepsilon$.
The proof that $\prob(\GeN^c)<\varepsilon$ is quite
technical and moved to the appendix.
Since $\varepsilon>0$ is arbitrary,
the above statements imply that
    $$
    \lim_{N\to \infty}\prob(H_{\sN}=3)=1.
    $$

Let $\varepsilon>0$ be fixed. The event
$\GeN$ is defined by
    $$
    \GeN =\{D_1^{\smallsup{N}}\le b_{\sD,\eps},\,
    D_2^{\smallsup{N}}\le b_{\sD,\eps},\, L^{\smallsup{N}}_{\sN}\ge
    L_{\varepsilon}(N)N^{1+\alpha(2-\tau)}\},
    $$
where, as before,  $b_{\sD,\eps} =\min\{k: 1-F(k)<\eps/8\}$ and
$L_{\varepsilon}(N)$ is some function satisfying for each
$\delta>0$ that
    $$
    \liminf_{N\to \infty} N^{\delta}L_{\eps}(N)>1.
    $$

\begin{lemma}
\label{technlemmatheorem1}
For each $\eps>0$,  there exists $N_{\eps}$, such that, for all $N\geq N_{\eps}$,
    \begin{equation}
    \label{afschattingGeN}
    \prob(\GeN^c)< \eps.
    \end{equation}
\end{lemma}
As explained above, the proof of this lemma is rather technical
and can be found in the appendix.

\medskip
\noindent
{\bf Proof of Theorem \ref{thm-tau(1,2)-cond}, part (i).}
We start with the  proof that $\prob\left(\{H_{\sN}=1\}\cap G_{\varepsilon, \sN}\right)<\varepsilon$, for
sufficiently large $N$. The event $\{H_{\sN}=1\}$ occurs {\it iff} at
least one stub of node $1$ connects to a stub of node $2$.
As before we denote for
$j\le D_1^{\smallsup{N}}$ by $\{[1.j]\to[2]\}$ the event
that $j$-th stub of node $1$ attaches to a stub of node $2$. Then
\begin{equation}
\label{hopcount=1}
\begin{array}{rl}
\prob\left(\{H_{\sN}=1\}\cap G_{\varepsilon, \sN}\right) \le& \expec\left[
\sumu\limits_{j=1}^{D_1^{\smallsup{N}}}
\prob_{\sN}\left([1.j]\to[2],\,G_{\varepsilon, \sN}\right)\right]\\[25pt]
\le&\expec\left[ \sumu\limits_{j=1}^{D_1^{\smallsup{N}}}
\frac{D_2^{\smallsup{N}}}{L^{\smallsup{N}}_{\sN}-1}\indic\{G_{\varepsilon, \sN}\}\right]
\le\frac{b_{\sD,\eps}^2}{N-1} <\varepsilon,
\end{array}
\end{equation}
for large enough $N$, since $L^{\smallsup{N}}_{\sN}\ge N$.

We will now bound $\prob\left(\{H_{\sN}=2\}\cap G_{\varepsilon, \sN}\right)$.
Note that the event $\{H_{\sN}=2\}$ occurs {\it iff} there exists node
$k\ge3$ and two stubs of node $k$ such that the first one connects
to a stub of node $1$ and the second one connects to a stub of
node $2$. For $k\ge3$ such that $D_k^{\smallsup{N}}\ge2$, $i,j\le
D_k^{\smallsup{N}}$, $i\ne j$, we denote by
$\{[k.i]\to[1],\,[k.j]\to[2]\}$ the event that the $i$-th stub of
node $k$ connects to a stub of node $1$ and the $j$-th stub of
node $k$ connects to a stub of node $2$. Then,
    \begin{equation}
    \label{hopcount=2}
    \begin{array}{rl}
    \prob\left(\{H_{\sN}=2\}\cap\GeN\right) \le& \expec\left[
    \sumu\limits_{k=3}^N\sumu\limits_{i\ne j}^{D_k^{\smallsup{N}}}
    \prob_{\sN}\left(\{[k.i]\to[1],\,[k.j]\to[2]\}\cap \GeN\right)\right]\\[25pt]
    \le&\expec\left[\sumu\limits_{k=3}^{N} \sumu\limits_{i\ne
    j}^{D_k^{\smallsup{N}}}
    \frac{D_1^{\smallsup{N}}}{L^{\smallsup{N}}_{\sN}-1}
    \frac{D_2^{\smallsup{N}}}{L^{\smallsup{N}}_{\sN}-3}
    \indic\{G_{\varepsilon, \sN}\}\right]\\
    \le&N \expec\left((D^{\smallsup{N}})^2\right)
    \frac{b_{\sD,\eps}^2}{(L_{\varepsilon}(N) N^{1+\alpha(2-\tau)}-3)^2}.
    \end{array}
    \end{equation}
Observe that
    $$
    \begin{array}{rl}
    \expec[(D^{\smallsup{N}})^2]
    &=\sum_{n< N^{\alpha}}(2n-1)\prob(D_1\geq n|D_1< N^{\alpha})\\
    &=\frac1{\prob(D_1<N^{\alpha})} \sum_{n< N^{\alpha}} (2n-1)L(n-1)(n-1)^{1-\tau}
    =L_4(N)N^{\alpha(3-\tau)},
    \end{array}
    $$
for some slowly varying function $L_4$. Substitution of this upper
bound in the right-hand side of~(\ref{hopcount=2})
shows that for $N$ large enough,
    $$
    \prob\left(\{H_{\sN}=2\}\cap \GeN\right)\leq
    N L_4(N)N^{\alpha(3-\tau)} \frac{b_{\sD,\eps}^2}{(L_{\varepsilon}(N)
    N^{1+\alpha(2-\tau)}-3)^2} <\varepsilon,
    $$
because $\alpha(3-\tau)<2(1+\alpha(2-\tau))$ when $\alpha<1/(1-\tau)$.

We will complete the proof by showing that, for sufficiently large $N$,
    \begin{equation}
    \label{bovengrens3}
        \prob(\{H_{\sN} \le 3\}\cap \GeN)\ge 1-\eps.
    \end{equation}
Let $\beta=(1+\alpha(4-\tau))/4$. Since $1/\tau<\alpha<1/(\tau-1)$,
we have
    \begin{equation}
    \label{ongebetalph}
    \frac{1+\alpha(2-\tau)}{2}<\beta<\alpha.
    \end{equation}
In this section we will call a node $k$, for $1\le k \le N,$ a {\it giant node}
if its degree $D_k^{\smallsup{N}}$ satisfies
    \begin{equation}
    \label{defgiant}
    N^{\beta}<D_k^{\smallsup{N}}< N^{\alpha}.
    \end{equation}
We will show below that with probability close to
$1$ at least one of the stubs of node~$1$ and at least one of the stubs of node~$2$
are connected to stubs of giant nodes, and that any two giant nodes have mutual
graph distance $1$. This
implies that the hopcount between the nodes 1 and 2 is at most
$3$.

The non-giant nodes, i.e. nodes with degree less than or equal to
$N^{\beta}$, are called {\em normal nodes} . First we will
show that the total degree of the normal nodes is negligible with
respect to $L^{\smallsup{N}}_{\sN}$. The mean degree of a normal node is
    $$
    \begin{array}{rl}
    \expec[D^{\smallsup{N}}\indic\{D\leq N^{\beta}\}]
    =&\sumu\limits_{n=1}^{\lfloor N^{\beta}\rfloor} \prob(D\geq n|D< N^{\alpha})\\
    =&\frac1{\prob(D< N^{\alpha})} \sumu\limits_{n=1}^{\lfloor
    N^{\beta}\rfloor} L(n-1)(n-1)^{1-\tau}
    =L_5(N)N^{\beta(2-\tau)},
    \end{array}
    $$
for some slowly varying function $L_5$. Thus, by the Markov
inequality,
    $$
    \prob\left(\sumu\limits_{i=1}^ND_i^{\smallsup{N}}\indic\{D_i\leq N^{\beta}\}
    \ge \frac{3}{\varepsilon}L_5(N)N^{1+\beta(2-\tau)}\right) \le
    \frac{\varepsilon}{3},
    $$
so that, with probability at least
$1-\varepsilon/3$, the fraction of the contribution from normal
nodes on $\GeN$, for sufficiently large $N$, is at
most
    $$
    \frac{
    \frac{3}{\varepsilon}
    L_5(N)N^{1+\beta(2-\tau)}}
    {L_{\varepsilon}(N)N^{1+\alpha(2-\tau)}}
    =\frac{3L_5(N)}
    {\eps L_{\eps}(N)}
    N^{2(\beta-\alpha)}.
    $$
Since $\beta<\alpha$, for
$\tau\in(1,2)$ the above fraction tends to $0$, as $N\to\infty$.
Thus the total contribution of the normal nodes is negligible with
respect to $L^{\smallsup{N}}_{\sN}$ on $\GeN$. This
implies that for sufficiently large $N$, with probability at most
$1-2\varepsilon/3,$ both nodes $1$ and $2$ are at graph distance $1$ from
some giant node on $\GeN$. It remains to show that any two giant nodes
have graph distance $1$ with
probability at least $1-\varepsilon/3$. For this we need an upper bound
of $L^{\smallsup{N}}_{\sN}$. Similarly as above we obtain
    $$
    \expec[D_1^{\smallsup{N}}]=L_6(N)N^{\alpha(2-\tau)},
    $$
for some slowly varying function $L_6$. Hence from the Markov inequality and since\\
$L^{\smallsup{N}}_{\sN}=D^{\smallsup{N}}_1+\dots+D^{\smallsup{N}}_{\sN}$,
$$
\prob\left(L^{\smallsup{N}}_{\sN}
>\frac{6}{\varepsilon}L_6(N)N^{1+\alpha(2-\tau)}\right)
\le
\frac{\eps N L_6(N)N^{\alpha(2-\tau)}}
{6L_6(N)N^{1+\alpha(2-\tau)}}\le\frac{\varepsilon}{6}.
$$
Since we have at most $N(N-1)<N^2$
pairs of giant nodes, the probability that two of them, say $g_1$ and
$g_2$, have graph distance greater than $1$ is at most (compare (4.5) of \cite{HHV03}),
\begin{eqnarray*}
&&\expec\left(N^2\produ\limits_{j=0}^{\lfloor D_{g_1}/2\rfloor-1}
\left(1-\frac{D_{g_2}}{L^{\smallsup{N}}_{\sN}-2i-1}\right)
\indic\{L^{\smallsup{N}}_{\sN}\le\frac{6}{\varepsilon}L_6(N)N^{1+\alpha(2-\tau)}\}\right)\\
&&\qquad+\prob\left((L^{\smallsup{N}}_{\sN}>\frac{6}{\varepsilon}L_6(N)N^{1+\alpha(2-\tau)}\right)
\le N^2\left(1-\frac{N^{\lfloor \beta\rfloor}}
{\frac{6}{\varepsilon}L_6(N)N^{1+\alpha(2-\tau)}}\right)^{\frac{1}{2}N^{\beta}}
+\frac{\varepsilon}{6}\\
&&\qquad=N^2
\left(e^{-1}+o(1)\right)
^{\eps N^{\frac
{2\beta-(1+\alpha(2-\tau))}
{(12L_6(N))}}}+
\frac{\varepsilon}{6}
\le\frac{\varepsilon}{3},
\end{eqnarray*}
for large $N$, since the exponent grows faster than any power of $N$.
Thus, with probability at least $1-\frac{\varepsilon}{3}$, all giant
nodes are on graph distance $1$, and thus form a complete graph.
This completes the proof.
\qed

\section{Simulation and conclusions}
\label{sec-sim}
To illustrate Theorem \ref{thm-tau(1,2)-uncond} and Theorem
\ref{thm-tau(1,2)-cond}, we have simulated our
random graph with degree distribution $D=\lceil
U^{-\frac{1}{\tau-1}} \rceil$, where $U$ is uniformly distributed
over $(0,1)$ and where for $x\in \mathbb R$, $\lceil x \rceil$ is
the smallest integer greater than or equal to $x$. Thus,
\begin{eqnarray*}
1-F(k) =\prob(U^{-\frac{1}{\tau-1}}>k)
=k^{1-\tau},\quad k=1,2,3,\ldots
\end{eqnarray*}

\renewcommand{\epsfsize}[2]{0.6#1}
\begin{figure}[t]
\begin{center}
\epsfbox[40 60 556 449]{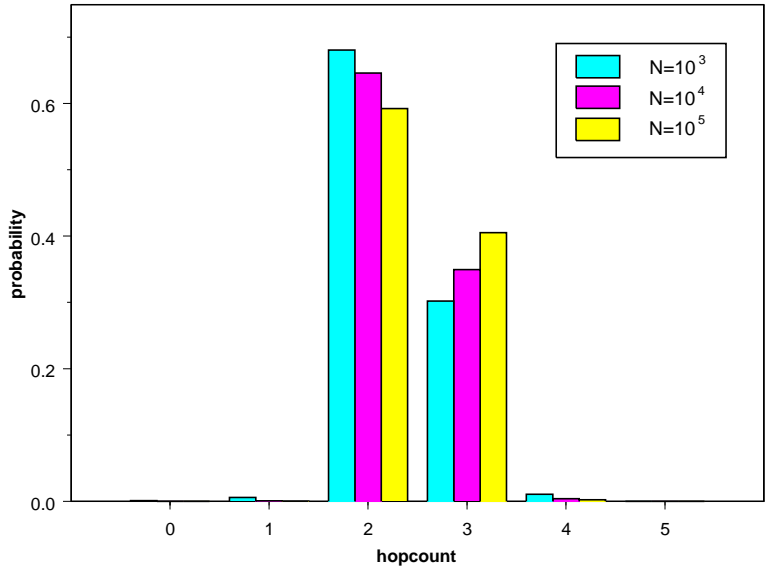}\caption{Empirical probability mass function
of the hopcount for $\tau=1.8$ and $N=10^3,10^4,10^5$, for the
unconditioned degrees.}
\label{histogram2}
\end{center}
\end{figure}

In Figure \ref {histogram2}, we have simulated the graph distance or hopcount
with $\tau=1.8$ and the values of
$N=10^3,10^4,10^5$. The histogram is in accordance with
Theorem \ref{thm-tau(1,2)-uncond}: for increasing values of $N$ we see that the
probability mass is divided over the values $H_{\sN}=2$ and $H_{\sN}=3$,
where the probability $\prob(H_N=2)$ converges.

\begin{figure}[t]
\begin{center}
\epsfbox[40 60 556 449]{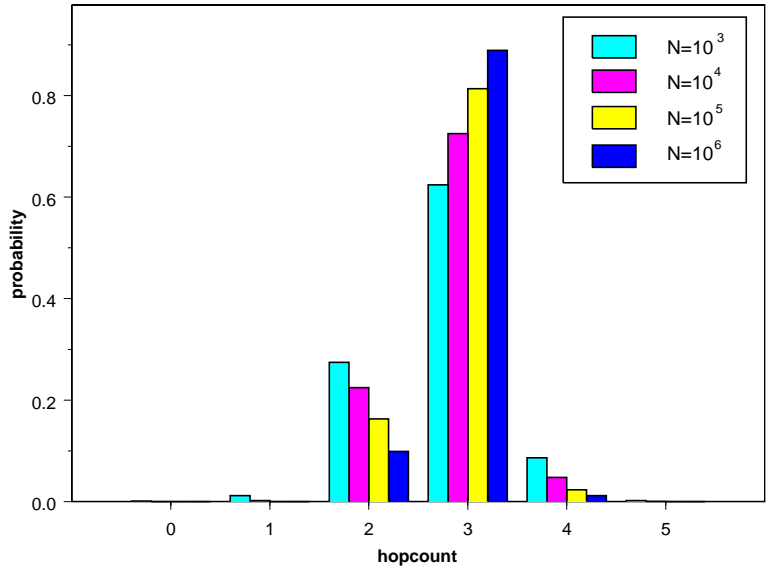}\caption{Empirical probability mass function
of the hopcount for $\tau=1.8$ and $N=10^3,10^4,10^5,10^6,$
where the degrees are conditioned to be less than $N$.}
\label{histogram1}
\end{center}
\end{figure}

As an illustration of Theorem \ref{thm-tau(1,2)-cond}, we again take $\tau=1.8$,
but now conditioned the degrees to be less than $N$, so that $\alpha=1$.
Since in this case $(\tau-1)^{-1}=\frac54$, we expect from
Theorem \ref{thm-tau(1,2)-cond} case (i), that in the limit the hopcount will concentrate
on the value $H_{\sN}=3$. This is indeed the case as is shown in Figure \ref{histogram1}

Our results give convincing asymptotics for
the hopcount when the mean degree is infinite, using extreme value theory.
Some interesting problems remain, such as
\begin{itemize}
\item[(i)]
Can one compute the exact value of $p_{\tau}$?

\item[(ii)]
What is the limit behavior of the hopcount when we condition the degrees on being
less than $N^{\alpha}$, with $\alpha<1/\tau$?

\end{itemize}

\renewcommand{\epsfsize}[2]{0.6#1}

\section*{Acknowledgement}
The work of RvdH and DZ was supported in part by Netherlands
Organisation for Scientific Research (NWO).



\renewcommand{\thesection}{\Alph{section}}
\setcounter{section}{0}

\numberwithin{equation}{subsection}
\numberwithin{theorem}{subsection}
\section{Appendix.}
\subsection{Proof of Lemma \ref{technlemma1}}
 In this section we restate Lemma
\ref{technlemma1} and then give a proof.

    \begin{lemma}
    \label{technlemma2}
    For each $\eps>0$, there exists $N_{\eps}$ such that
        \begin{equation}
        \label{afschattingAeN2}
        \prob(\AeN^c)< \eps, \quad N\geq N_{\eps}.
        \end{equation}
    \end{lemma}

\noindent{\bf Proof.}
We define the event $\EeN$ by
    \begin{equation}
    \label{dnz7}
    \begin{array}{rll}
    \EeN=
    &\left\{\sum\limits_{n=0}^{k_\varepsilon}D_{(N-n)}
    \le\frac{\varepsilon}{8b_{\sD,\eps}}L_{\sN}\right\}\qquad&(a)\\
    &\cap\left\{D_{\smallsup{N-k_\varepsilon+1}}\ge a_{\xi_{k_\eps}}u_{\sN}\right\}\qquad&(b)\\
    &\cap\left\{L_{\sN}\le b_{\eta, \eps}u_{\sN}\right\},\qquad&(c)
    \end{array}
    \end{equation}
where $b_{\sD,\eps}$ is the $\eps$-quantile of $F$ used in the definition of
$\DeN$ and where $a_{\xi_{k_\eps}},b_{\eta, \eps}>0$ are defined by
    $$
    \prob\left(\xi_{k_\eps}<a_{\xi_{k_\eps}}\right) < \eps/24
    \quad \text{and} \quad
    \prob(\eta>b_{\eta, \eps})<\eps/24,
    $$
respectively. Observe that $a_{\xi_{k_\eps}}$ is a lower quantile of $\xi_{k_\eps}$,
whereas $b_{\xi_{k_\eps}}$ defined in (\ref{defa4N}) and (\ref{defkeps})
is an upper quantile of $\xi_{k_\eps}$.
Furthermore, $b_{\eta, \eps}$ is an upper quantile of $\eta$,
whereas $a_{\eta, \eps}$ defined in (\ref{defa3N}) is a lower
quantile of $\eta$.
Since
    $$
    \AeN=\BeN\cap\CeN\cap\DeN,
    $$
(see (\ref{def-AeN}) and below the proof of Lemma  \ref{Ldnz1} for the definition of
$\AeN$, $\BeN$, $\CeN$ and $\DeN$). we have
    \begin{equation}
    \label{dnz67_4}
    \prob(\AeN^c)\le \prob(\BeN^c\cap\DeN\cap\EeN)
    +\prob(\CeN^c\cap\DeN\cap\EeN)
    +\prob(\DeN^c)+\prob(\EeN^c),
    \end{equation}
and in order to prove the lemma we should show that each of the four terms
on the right-hand side is at most $\eps/4$.

Since on $\DeN$, the nodes 1 and 2 each have at most $b_{\sD,\eps}$ stubs,
the first term satisfies
    $$
    \prob(\BeN^c\cap\DeN\cap\EeN) \le 2b_{\sD,\eps}\expec\left(
    \frac{1}{L_{\sN}}\sum\limits_{n=0}^{k_{\varepsilon}} D_{\smallsup{N-n}}
    \indic\{\EeN\}\right) \le\varepsilon/4,
    $$
due to point $(a)$ of $\EeN$. This bounds the first term of (\ref{dnz67_4}).

We turn to the second term of (\ref{dnz67_4}). Recall that $\CeN^c$ induces that no stubs of
at least two giant nodes are attached to one another.
Since we have at most $N^2$ pairs of giant
nodes $g_1$ and $g_2$, the items $(b)$, $(c)$ of $\EeN$ imply
    $$
    \begin{array}{rl}
    \prob(\CeN^c\cap\DeN \cap\EeN)\le
    &\expec\left(N^2\produ\limits_{j=0}^{\lfloor D_{g_1}/2\rfloor-1}
    \left(1-\frac{D_{g_2}}{L_{\sN}-2i-1}\right)\right)\\[10pt]
    \le&N^2\left(1-\frac{a_{\xi_{k_\eps}} u_{\sN}} {b_{\eta, \eps}
    u_{\sN}}\right)
    ^{a_{\xi_{k_\eps}}u_{\sN}/2}\\[10pt]
    =&N^2\left(e^{-1}+o(1)\right)
    ^{(a_{\xi_{k_\eps}})^2/(2b_{\eta, \eps}) u_{\sN}}\\
    =&N^2\left(e^{-1}+o(1)\right)
    ^{(a_{\xi_{k_\eps}})^2/(2b_{\eta, \eps})L^{'}(N)N^{\frac{1}{\tau-1}}}
    \le\frac{\varepsilon}{4},
    \end{array}
    $$
for large enough $N$,
since the exponent grows faster than any power of $N$.

The third term at the r.h.s. of~(\ref{dnz67_4}) is at most
$\varepsilon/4$, because
    $$
    \prob(\DeN^c)\le 2 \prob(D_1>b_{\sD,\eps})
    \le2\varepsilon/8=\eps/4.
    $$

It remains to estimate the last term
at the r.h.s. of~(\ref{dnz67_4}). Clearly,
\begin{equation}
\label{dnz8}
\begin{array}{rll}
\prob\left(\EeN^c\right)\le
&\prob\left(\sum\limits_{n=0}^{k_\varepsilon}D_{(N-n)}
>\frac{\varepsilon}{8b_{\sD,\eps}}L_{\sN}\right)\qquad&(a)\\
&+\prob\left(D_{\smallsup{N-k_\varepsilon+1}}<a_{\xi_{k_\eps}}u_{\sN}\right)\qquad&(b)\\
&+\prob\left(L_{\sN}>b_{\eta, \eps}u_{\sN}\right).\qquad&(c)
\end{array}
\end{equation}
We will consequently  show that
each term in the above expression is at most $\varepsilon/12$.
Let $a_{\eta, \eps}$ and $b_{\xi_{k_\eps}}>0$
be as in (\ref{defa3N}) and (\ref{defa4N}), then we can
decompose the first term at the r.h.s. of~(\ref{dnz8}) as
    \begin{eqnarray}
    \label{dnz9}
    \prob\left(\sum\limits_{n=0}^{k_\varepsilon}D_{\smallsup{N-n}}
    >\frac{\varepsilon}{8b_{\sD,\eps}}L_{\sN}\right)&\le&
    \prob\left(L_{\sN}<a_{\eta, \eps}u_{\sN}\right)+\prob\left(\sum\limits_{n=0}^{k_\varepsilon}D_{\smallsup{N-n}}
    >\frac{\varepsilon}{8b_{\sD,\eps}}a_{\eta, \eps}u_{\sN}\right)\\
    &\le&
    \prob\left(L_{\sN}<a_{\eta, \eps}u_{\sN}\right)+\prob\left(D_{\smallsup{N-k_\varepsilon+1}}>b_{\xi_{k_\eps}}u_{\sN}\right)\nn\\
    &&\qquad
    +\prob\left(\sum\limits_{i=1}^N D_i\indic\{D_i<b_{\xi_{k_\eps}}u_{\sN}\}
    >\frac{\varepsilon}{8b_{\sD,\eps}}a_{\eta, \eps}u_{\sN}\right).\nonumber
    \end{eqnarray}
From the Markov inequality,
    \begin{equation}
    \label{dnz14} \prob\left(\sum\limits_{i=1}^N
    D_i\indic\{D_i<b_{\xi_{k_\eps}}u_{\sN}\}
    >\frac{\varepsilon}{8b_{\sD,\eps}}a_{\eta, \eps}u_{\sN}\right)
    \le \frac {N\expec\left(D_1\indic\{D_1<b_{\xi_{k_\eps}}
    u_{\sN}\}\right)} {\frac{\varepsilon}{8b_{\sD,\eps}}a_{\eta, \eps}u_{\sN}}.
    \end{equation}
Since $1-F(x)$ varies regularly with exponent $\tau-1$ we have,
by~\cite[ Theorem 1(b), VIII.9, p.281]{fellerb},
    \begin{equation}
    \label{dnz15}
    \expec\left(D_1\indic\{D_1<b_{\xi_{k_\eps}}u_{\sN}\}\right)
    =\sumu\limits_{k=0}^{\lfloor b_{\xi_{k_\eps}}u_{\sN}\rfloor}
    \left(1-F(k)\right)\le 2(2-\tau)b_{\xi_{k_\eps}}u_{\sN}
    \left(1-F(b_{\xi_{k_\eps}}u_{\sN})\right),
    \end{equation}
for large enough $N$. Due to~(\ref{dnz11}), for large enough $N$,
we have also
    \begin{equation}
    \label{dnz18}
    N\left(1-F(u_{\sN})\right)\leq 2.
    \end{equation}
Substituting (\ref{dnz15}) and~(\ref{dnz18})
in~(\ref{dnz14}), we obtain
\begin{eqnarray}
    \label{dnz19}
    \prob\left(
    \sum_{i=1}^N
    D_i\indic\{D_i<b_{\xi_{k_\eps}} u_{\sN}\}
    >\frac{\varepsilon}{8b_{\sD,\eps}}a_{\eta, \eps}u_{\sN}\right)
    &\leq&
    \frac{2N(2-\tau)b_{\xi_{k_\eps}}u_{\sN}
    \left(1-F(b_{\xi_{k_\eps}}u_{\sN})\right)}{\frac{\varepsilon}{8b_{\sD,\eps}}u_{\sN}a_{\eta, \eps}}
   \nonumber\\
        &\le&\frac{4(2-\tau)b_{\xi_{k_\eps}}
    \left(1-F(b_{\xi_{k_\eps}}u_{\sN})\right)}
    {\frac{\varepsilon}{8b_{\sD,\eps}}a_{\eta, \eps}\left(1-F(u_{\sN})\right)},
\end{eqnarray}
for large enough $N$. Since $1-F(x)$ varies regularly with
exponent $\tau-1$,
    $$\lim_{N\to\infty}
    \frac {\left(1-F(b_{\xi_{k_\eps}}u_{\sN})\right)}
    {\left(1-F(u_{\sN})\right)}=
    \left(b_{\xi_{k_\eps}}\right)^{1-\tau}.
    $$
Hence the r.h.s. of~(\ref{dnz19}) is at most
    $$
    \frac {8(2-\tau)\left(b_{\xi_{k_\eps}}\right)^{2-\tau}}
    {\frac{\varepsilon}{8b_{\sD,\eps}}a_{\eta, \eps}} =
    \varepsilon/36,
    $$
for sufficiently large $N$, by definition of $b_{\xi_{k_\eps}}$ in (\ref{defa4N}).
We now show that the second term of~(\ref{dnz9}) is at most $\varepsilon/36$.
Since $D_{\smallsup{N-k_\varepsilon+1}}/u_{\sN}$ converges in distribution to $\xi_{k_\eps}$,
we find from (\ref{defkeps}),
    $$
    \prob\left(D_{\smallsup{N-k_\varepsilon+1}}>b_{\xi_{k_\eps}}u_{\sN}\right)
    \le\prob\left(\xi_{k_\varepsilon}>b_{\xi_{k_\eps}}/2\right)
    +\eps/72
    \le\varepsilon/36,
    $$
for large enough $N$. Similarly, by definition of $a_{\eta, \eps}$, in (\ref{defa3N}),  we have
    $$
    \prob\left(L_{\sN}<a_{\eta, \eps}u_{\sN}\right)\le \eps/36.
    $$
Thus, the term~(\ref{dnz8})(a) is at most $\varepsilon/12$.

The upper bound for~(\ref{dnz8})(b), i.e., the bound
    $$
    \prob\left(D_{\smallsup{N-k_\varepsilon+1}}<a_{\xi_{k_\eps}}u_{\sN}\right)<\eps/12,
    $$
is an easy consequence of the distributional convergence
of  $D_{\smallsup{N-k_\eps+1}}/u_{\sN}$ to $\xi_{k_\eps}$
and the definition of $a_{\xi_{k_\eps}}$.
Similarly, we obtain the upper bound for the term in~(\ref{dnz8})(c), i.e.,
    $$
    \prob\left(L_{\sN}>b_{\eta, \eps}u_{\sN}\right)<\eps/12,
    $$
from the distributional convergence of $L_{\sN}/u_{\sN}$ to $\eta$ and the definition of
$b_{\eta, \eps}$.

Thus we have shown that $\prob(\EeN^c)<\varepsilon/4$.
This completes the proof of Lemma \ref{technlemma1}.
\qed

\subsection{Proof of Lemma \ref{technlemmatheorem1}}

In this section we restate Lemma \ref{technlemmatheorem1} and give a proof.
    \begin{lemma}
    \label{technlemmatheorem2}
    For each $\eps>0$, there exists $N_{\eps}$ such that for all $N\geq N_{\eps}$,
        \begin{equation}
        \label{afschattingGeN2}
        \prob(\GeN^c)< \eps.
        \end{equation}
    \end{lemma}

\noindent{\bf Proof.}
Clearly,
    \begin{equation}
    \label{dnz1}
    \prob(\GeN^c)<2\prob\left(D_1^{\smallsup{N}}>b_{\sD,\eps}\right)
    +\prob\left(L^{\smallsup{N}}_{\sN}<L_{\varepsilon}(N)N^{1+\alpha(2-\tau)}\right).
    \end{equation}
From the definition of $b_{\sD,\eps}\in\NN$, and the string of inequalities,
    $$
    \prob(D_1^{\smallsup{N}}>b_{\sD,\eps})
    \le\prob(D_1>b_{\sD,\eps})<\varepsilon/4,
    $$
we obtain that the first term on the right-hand side of ~(\ref{dnz1}) is at most  $\eps/2$.
We show now that the second term on the right-hand side of~(\ref{dnz1}) is at most
$\eps/2$.
Since
    $$
    D^{\smallsup{N}}_i=\sum_{j<N^{\alpha}}
    \indic\{D^{\smallsup{N}}_i\ge j\},
    $$
we obtain, after interchanging the order of
the two summations involved, that
    \begin{equation}
    \label{dnz4}
    L^{\smallsup{N}}_{\sN}=\sum_{i=1}^ND^{\smallsup{N}}_i
    =\sum_{j=1}^{\lfloor N^{\alpha}\rfloor-1}\sumu\limits_{i=1}^N\indic\{D^{\smallsup{N}}_i\ge j\}.
    \end{equation}
We would like to approximate
    $$
    \sumu\limits_{i=1}^N\indic\{D^{\smallsup{N}}_i\ge j\}\qquad \mbox{ by }\qquad
    N\prob\left(D^{\smallsup{N}}_1\ge j\right).$$
For this we will use
(\cite[Theorem 1.7(i), p.~14]{bol01}) which states
that for a binomial random variable $X$ with parameters $N$ and
$p\leq 1/2$, and $a$ such that
    \begin{equation}
    \label{arest}
    a(1-p)\geq 12 \qquad \mbox{and} \qquad 0<\frac{a}{Np}\leq \frac{1}{12},
    \end{equation}
we have
    \begin{equation}
    \prob(|X-Np|\geq a) \leq \sqrt{\frac{Np}{a^2}} e^{-\frac{a^2}{3Np}}.
    \label{binbd}
    \end{equation}
Observe that $p_{\sN}(j)=\prob\left(D^{\smallsup{N}}_1\ge j\right)$
is  non-increasing in $j\geq 1$ and that
$\lim\limits_{N\to\infty}p_{\sN}(N-1)=0$.
We will first complete the proof under an extra assumption on $\{p_{\sN}(j)\}_{j\geq 1},$
and then prove that the assumption indeed holds.
The extra assumption reads that for some $1\le\kappa<k_{\sN}<N^{\alpha}$, all $j\in[\kappa,k_{\sN}]$
and $N$ large enough, we have
    \begin{equation}
    \label{dnz2}
    \begin{array}{rl}
    (a)&\qquad p_{\sN}(j)<1/2,\\
    (b)&\qquad N p_{\sN}(j)>288,\\
    (c)&\qquad N p_{\sN}(j)/432-\alpha\log N\to\infty,\mbox{ as }N\to\infty,\\
    (d)&\qquad p_{\sN}(j)>\frac{1}{2}\prob\left(D_1\ge j\right).
    \end{array}
    \end{equation}
Then, with $a=Np_{\sN}(j)/12$, we have that
    \begin{eqnarray*}
    \prob\left(\sum_{i=1}^{N
    }\indic\{D^{\smallsup{N}}_i <
    j\} < Np_{\sN}(j)/2)\right) &\leq&
    \prob
    (
    |
    \sum_{i=1}^{N
    }\indic\{D^{\smallsup{N}}_i <
    j\}-Np_{\sN}(j))
    |>a/2))\\
    &\leq& 12(Np_{\sN}(j))^{-1/2}e^{-Np_{\sN}(j)/432},
    \end{eqnarray*}
because from (a) and (b) of (\ref{dnz2}), we have
$a(1-p)=(Np_{\sN}(j)/12)(1-p_{\sN}(j))\ge Np_{\sN}(j)/24>288/24=12$ and $a/Np=1/12$.
Hence, we find
    \begin{eqnarray}
    \label{dnz3}
    &&\prob\left(\bigcap_{j=\kappa}^{k_{\sN}}\{
    \sum_{i=1}^{N}
    \indic\{D^{\smallsup{N}}_i \geq j\} \geq Np_{\sN}(j)/2\}\right)
    =1-\prob\left(\bigcup_{j=\kappa}^{k_{\sN}}\{
    \sum_{i=1}^N
    \indic\{D^{\smallsup{N}}_i < j\} \geq Np_{\sN}(j)/2\}\right)\nonumber\\
    &&\quad
    \geq 1-\sum_{j=\kappa}^{k_{\sN}}\prob\left(\{
    \sum_{i=1}^N
    \indic\{D^{\smallsup{N}}_i < j\} \geq Np_{\sN}(j)/2\}\right)\nonumber\\
    &&\quad
    \geq 1-\sum_{j=\kappa}^{k_{\sN}} 12(N p_{\sN}(j))^{-1/2}\exp(-Np_{\sN}(j)/432)
    \geq 1- k_{\sN}\exp(-Np_{\sN}(k_{\sN})/432)\nonumber\\
    &&\quad
    \geq 1-\exp\{\alpha\log(N)-N p_{\sN}(k_{\sN})/432\}\geq 1- \eps/2,
    \end{eqnarray}
for large $N$,
due to~(\ref{dnz2}$(b)$) and~(\ref{dnz2}$(c)$). Given~(\ref{dnz3})
and~(\ref{dnz2},$(d)$), we have
    \begin{equation}
    \label{dnz5}
    L_{\sN}^{\smallsup{N}}\geq\sum_{j=\kappa}^{k_{\sN}}
    \sum_{i=1}^N
    \indic\{D^{\smallsup{N}}_i\ge j\}
    \ge\frac{N}{2}\sum_{j=\kappa}^{k_{\sN}}p_{\sN}(j)
    \ge\frac{N}{4}\sum_{j=\kappa}^{k_{\sN}}\prob(D\ge j)=\frac{N}{4}
    \sum_{j=\kappa-1}^{k_{\sN}-1} j^{1-\tau}L(j),
    \end{equation}
with probability at least $1-\eps/2$.

We will call a function $f(N)$ {\em slow}, if for each
$\delta>0$, $f(N)>N^{-\delta}$, for large enough $N$. Observe that
for any $a>0$ and slow $f$
    $$
    \sum_{j=1}^kf(j)j^a=f_1(k)k^{a+1},
    $$
for some slow function $f_1$.
We further assume that we can take $k_{\sN}=L_1(N)N^{\alpha}$ for some slow
function $L_1(N)$, then the r.h.s. of~(\ref{dnz5}) is
$L_{\varepsilon}(N)N^{1+\alpha(2-\tau)}$, for some slow function
$L_{\varepsilon}(N)$. This would imply that the second term
of~(\ref{dnz1}) is at most $\varepsilon/2$ and thus completes the
proof of Lemma \ref{technlemmatheorem2} subject to
(\ref{dnz2}) and the fact that $k_{\sN}=L_1(N)N^{\alpha}$ for some slow
function $L_1(N)$.

We now specify
$1\le\kappa<N^\alpha$, check
points~(\ref{dnz2}$(a)$)-~(\ref{dnz2}$(d)$) and demonstrate that
we can take $k_{\sN}=L_1(N)N^{\alpha}$, for some slow
function $L_1(N)$. The only restriction
on $\kappa$ is $(\ref{dnz2}$(a)$)$. Just take $\kappa$ large
enough such that $p_{\sN}(\kappa)<1/2$. Then, since $j\mapsto p_{\sN}(j)$
is non-increasing, we obtain that
$p_{\sN}(j)<1/2$ for all $j\ge\kappa$.
Therefore, (\ref{dnz2}$(a)$) is satisfied.
This introduces $\kappa$ and proves (\ref{dnz2}$(a)$).

We next define $k_{\sN}$. Choose
    $$
    k_{\sN}=\max\limits_{k}\{
    \prob(D_1>k)>2(1-F(\lfloor N^{\alpha}\rfloor-1))\}.
    $$
This definition gives us point $(d)$ of~(\ref{dnz2}). Indeed for
any $j\le k_{\sN}$,
    \begin{eqnarray}
    \label{dnz6}
    &&p_{\sN}(j)=\prob\left(D^{\smallsup{N}}_1\ge j\right)
    =\frac{\prob\left(D_1\ge j\right)-\prob\left(D_1\ge N^\alpha\right)}
    {1-\prob(D_1\ge N^\alpha)}\nonumber\\
    &&\qquad>\prob\left(D_1\ge j\right)-\prob\left(D_1\ge N^\alpha\right)
    \ge\frac{1}{2}\prob\left(D_1\ge j\right).
    \end{eqnarray}
Before we check the other items of~(\ref{dnz2}),
we prove that $k_{\sN}=L_1(N)N^{\alpha}$, for some slow function
$L_1(N)$. We argue by contradiction. Suppose there exists
$\delta>0$ such that $k_{\sN}<N^{(1-\delta)\alpha}$ for all $N$
sufficiently large. Then, by
definition of $k_{\sN}$, we have
    $$
    \prob(D_1>k_{\sN}+1)\le
    2\prob(D_1>\lfloor N^{\alpha}\rfloor)=L(N)N^{\alpha(1-\tau)}.
    $$
However, when we compute
probability on the l.h.s., then we obtain that
    $$
    \prob(D_1>k_{\sN}+1)\ge
    L_3(N)N^{-(\tau-1)(1-\delta)\alpha},
    $$
for some slow function $L_3(N)$.
This gives a contradiction. Hence, for any $\delta>0$,
$k_{\sN}>N^{\alpha(1-\delta)}$, for large enough $N$. This is equivalent to
the statement that $k_{\sN}=L_1(N)N^{\alpha}$ for some slow function $L_1(N)$.

Since $L_1(N)$ is slow, $(\ref{dnz2}$(c)$)$ follows from
$(\ref{dnz2}$(d)$)$, since
    \begin{equation}
    \label{dnz20}
    \frac{N}{2}\prob\left(D_1\ge
    k_{\sN}\right)=\frac{N}{2}L(k_{\sN}-1)(k_{\sN}-1)^{1-\tau}=L_2(N)N^{1-\alpha(\tau-1)},
    \end{equation}
for some slow function $L_2(N)$. Since $\alpha<1/(\tau-1)$, the right-hand side of
(\ref{dnz20}) increases as a positive power of $N$, and hence exceeds $\alpha \log(N)$,
eventually for sufficiently large $N$.
Finally, (\ref{dnz2}$(b)$) follows from~(\ref{dnz2}$(c)$). Thus, we have proved that
$\kappa,k_{\sN}$ defined above satisfy the assumptions in (\ref{dnz2}) and that
$k_{\sN}=L_1(N)N^{\alpha}$ for some slow
function $L_1(N)$. Therefore, the second term in~(\ref{dnz1}) is at most $\varepsilon/2$,
and we conclude that $\prob(\GeN^c)<\varepsilon$.
\qed

\end{document}